\documentclass[reqno,11pt] {amsart}

\newtheorem{theorem}{Theorem}[section]
\newtheorem{lemma}[theorem]{Lemma}

\theoremstyle{definition}

\theoremstyle{remark}
\newtheorem{remark}[theorem]{Remark}

\newcommand{\mysection}[1]{\section{#1}
\setcounter{equation}{0}}

\newcommand{\bR}{\mathbb R}

\newcommand\cF{\mathcal{F}}

\renewcommand{\epsilon}{\varepsilon}

\begin{document}
\title[On  Peng's central limit theorem]
 {On Shige Peng's central limit theorem}

\author[N.V. Krylov]{N.V. Krylov}

\address[N.V. Krylov]
{127 Vincent Hall, University of Minnesota, Minneapolis, MN
55455, USA}
 
\email{nkrylov@umn.edu}



\begin{abstract} 
We give  error estimates in 
Peng's central limit theorem
for not necessarily nondegenerate
case. The exposition uses the language of the
classical probability theory
instead of the language of the theory
of sublinear expectations.
We only consider the one-dimensional case.
The higher dimensional extension is straightforward.

\end{abstract}

\maketitle

\mysection{Introduction and the main result}
                                  \label{section 6.28.2}

Let $\Theta$ be a collection of real-valued
random variables, given perhaps on different
probability spaces such that
$$
 E\xi=0\quad\forall\xi\in\Theta,\quad
 \bar \sigma:= \sup_{\xi\in\Theta}\big[E(\xi^{2})\big]^{1/2}<\infty.
$$
Denote
$$
 \underline{\sigma}=\inf_{\xi\in\Theta}\big[E(\xi^{2})\big]^{1/2}.
$$
Fix a number $\beta\in(0,1]$ and suppose that
\begin{equation}
                                                  \label{6.27.1}
M_{\beta}:=\sup_{\xi\in\Theta} E(|\xi|^{2+\beta})<\infty.
\end{equation}

Next, let $(\Omega,\cF,P)$ be a complete probability space 
with an increasing
 filtration of $\sigma$-fields $\cF_{t}\subset\cF$, $t\geq0$,
each of which contains all $P$-null sets. Suppose that
on $\Omega$ we are given a one-dimensional standard
Wiener process $w_{t},t\geq0$, which is a Wiener process
with respect to $\{\cF_{t},t\geq0\}$. 
Let $\Sigma$ be the set
of   real-valued, progressively
measurable processes $\sigma_{\!\cdot}=\{\sigma_{s},s\in[0,1]\}$ satisfying
$\bar\sigma\geq \sigma_{s}\geq\underline{\sigma}$
for all $s$ and $\omega$. For $0\leq t\leq s\leq 1,x\in\bR^{d}$,
and $\sigma_{\!\cdot}\in\Sigma$, introduce
$$
x^{\sigma_{\!\cdot}}_{s}(t,x)=x+\int_{t}^{s}\sigma_{r}\,dw_{r}.
$$

Fix   
 a real-valued function $\phi$
on $\bR$ such that, for any $x,y\in\bR$,
\begin{equation}
                                                    \label{6.22.3}
|\phi(x)-\phi(y)|\leq |x-y|^{\beta}.
\end{equation}
Then, for $t\in[0,1],x\in\bR^{d}$, introduce 
\begin{equation}
                                                    \label{6.22.1}
v(t,x)=\sup_{\sigma_{\!\cdot}\in\Sigma}
 E\phi\big(x^{\sigma_{\!\cdot}}_{1}(t,x)\big).
\end{equation}

Finally, for   $n=1,2,...$, $k=0,...,n-1$,
  define recursively, $v_{n}(1,x)=\phi(x)$,
\begin{equation}
                                                    \label{6.22.4}
v_{n}( k/n,x)=
\sup_{\xi\in\Theta} Ev_{n}\big( (k+1)/n,x+\xi/\sqrt{n}\big).
\end{equation}

The first goal of this note is to prove the following theorem.

\begin{theorem}
                                         \label{theorem 6.22.1}
We have
\begin{equation}
                                                    \label{6.22.2}
|v(0,0)-v_{n}(0,0)|
\leq Nn^ {-\beta^{2}/(4+2\beta)} ,
\end{equation}
where $N$ is a constant depending only on $\bar\sigma$,
$M_{\beta}$, and $\beta$.
\end{theorem}

\begin{remark}
                                       \label{remark 6.27.1}
S. Peng (\cite{Pe_08}) proved that $v_{n}(0,0)\to v(0,0)$
in multidimensional setting when $\xi$ can contain drift terms
of order $1/\sqrt{n}$. He used his theory of sublinear expectations.
In the one-dimensional case with $E\xi=0$
our setting is {\em equivalent\/} to his.
His arguments rely on
 the theory of fully nonlinear parabolic equations
and assumption \eqref{6.27.1}, albeit, never stated explicitly.

Y. Song (\cite{So_17}) gave an estimate for the rate of convergence
similar to \eqref{6.22.2}, however, assuming
 that $\underline{\sigma}>0$ and with and exponent
which apparently should vanish as $\underline{\sigma}
\downarrow0$. Y. Song also assumes \eqref{6.27.1}
and uses the theory of fully nonlinear parabolic equations
and the theory of sublinear expectations.

Our approach is different and is based on an idea
from \cite{Kr_97}. The case that  $\underline{\sigma}=0$
is not excluded.

\end{remark}

We prove Theorem \ref{theorem 6.22.1} 
in Section \ref{section 6.26.1}
and in the rest of the article by $N$ we denote
generally different constants depending only on 
$\bar\sigma$,
$M_{\beta}$, and $\beta$.

In Section \ref{section 6.28.1} we present
a better rate of convergence of $v_{n}$ to $v$
under more restrictive conditions.

\mysection{Auxiliary estimates}

Here are a few properties of $v$ (see Sections 3.1
and 4.1 in \cite{Kr_77}).

\begin{lemma}
                                   \label{lemma 6.22.1}
(a) Bellman's principle holds:
for any $0\leq s\leq t\leq 1$
and $x\in\bR$,
$$
v(s,x)=\sup_{\sigma_{\cdot}\in\Sigma}
 Ev\big(t,x^{\sigma_{\!\cdot}}_{t}(s,x)\big).
$$

(b) For any $t\in[0,1]$, $x,y\in\bR$,
$$
|v(t,x)-v(t,y)|\leq |x-y|^{\beta}.
$$

(c)  For any $s,t\in[0,1]$, $x \in\bR$,
$$
|v(t,x)-v(s,x)|\leq \bar{\sigma}^{\beta}|t-s|^{\beta/2}.
$$
\end{lemma}

Proof. Assertion (a) is proved in \cite{Kr_77}.
Assertion (b) is a trivial consequence of \eqref{6.22.3}
and the fact that $x^{\sigma_{\!\cdot}}_{1}(t,x)-
x^{\sigma_{\!\cdot}}_{1}(t,y)=x-y$. Assertion (c) is given
as an exercise in \cite{Kr_77}, the solution 
of which in our particular case is as follows.

We may assume that $s\leq t$. Then by (a) and (b)
$$
|v(t,x)-v(s,x)|\leq 
\sup_{\sigma_{\cdot}\in\Sigma}
 E\big|v\big(x^{\sigma_{\!\cdot}}_{t}(s,x)\big)-v(s,x)\big|
$$
$$
\leq \sup_{\sigma_{\cdot}\in\Sigma}
 E\Big|\int_{s}^{t}\sigma_{r}\,dw_{r}\Big|^{\beta}
\leq \sup_{\sigma_{\cdot}\in\Sigma}
\Big( E\Big|\int_{s}^{t}\sigma_{r}\,dw_{r}\Big|^{2}\Big)^{\beta/2}
=|t-s|^{\beta/2}\bar{\sigma}^{\beta}.
$$
The lemma is proved.

Define $T_{n}=\{k/n:k=0,...,n\}$.

\begin{lemma}
                                   \label{lemma 6.22.2}
 
(a) For any $t\in T_{n}$, $x,y\in\bR$,
$$
|v_{n}(t,x)-v_{n}(t,y)|\leq |x-y|^{\beta}.
$$

(b) For any $t,s\in T_{n}$, $x \in\bR$,
$$
|v_{n}(t,x)-v_{n}(s,x)|\leq \bar{\sigma}^{\beta}|t-s|^{\beta/2}.
$$
\end{lemma}

Proof. Assertion (a) is proved by backward
induction using
\eqref{6.22.4}. 

To prove (b)
we concentrate on $s\leq t$, $t>0$, $s,t\in T_{n}$,  $x \in\bR$,
and
first claim that, if for some  
$a,b\geq 0$,
for any $ y\in\bR$,    we have
$$
v_{n}(t,x)\leq v_{n}(t,y)+a|x-y|^{2}+b,
$$
then for any $s\in T_{n},s\leq t$ we have
\begin{equation}
                                                   \label{6.25.1}
v_{n}(s,x)\leq v_{n}(t,y)+a|x-y|^{2}+a\bar{\sigma}^{2}(t-s)+b.
\end{equation}

Indeed, \eqref{6.25.1} holds for $s=t$ and if it holds for
an  $s\in T_{n}$, $0<s\leq t$,  then
$$
v_{n}(s-1/n,x)=\sup_{\xi\in\Theta}E v_{n}(s,x+\xi/\sqrt{n})
\leq v_{n}(t,y)+a\sup_{\xi\in\Theta}E|x+\xi/\sqrt{n}-y|^{2}
$$
$$
+a\bar{\sigma}^{2}(t-s)+b
=v_{n}(t,y)+a|x-y|^{2}+a\bar{\sigma}^{2}(t-s+1/n)+b.
$$
This proves our claim.

Since for any $a\geq0$, $\varepsilon>0$,
$$
a^{\beta}\leq (\beta/2)\varepsilon a^{2}+(1-\beta/2)\varepsilon
^{-\beta/(2-\beta)},
$$
we have
$$
v_{n}(t,x)\leq v_{n}(t,y)+|x-y|^{\beta}
\leq v_{n}(t,y)+(\beta/2)\varepsilon|x-y|^{2}
+(1-\beta/2)\varepsilon
^{-\beta/(2-\beta)},
$$
$$
v_{n}(s,x)\leq v_{n}(t,y)+(\beta/2)\varepsilon
|x-y|^{2}+ (\beta/2)\varepsilon\bar{\sigma}^{2}(t-s)
+(1-\beta/2)\varepsilon
^{-\beta/(2-\beta)}.
$$
 
By taking here $y=x$ and computing the inf of the right-hand side
with respect to $\varepsilon>0$ we get
$$
v_{n}(s,x)\leq v_{n}(t,x)+\bar{\sigma}^{\beta}(t-s)^{\beta/2}.
$$
Similarly one proves that
$$
v_{n}(s,x)\geq v_{n}(t,x)-\bar{\sigma}^{\beta}(t-s)^{\beta/2}
$$
and this yields (b). The lemma is proved.

 Take a nonnegative $\zeta\in C^{\infty}_{0}(\bR^{2})$
with unit integral and support in $\{(t,x):-1<t<0,
|x|<1\}$ and for $\varepsilon\in(0,1)$ introduce
$\zeta_{\varepsilon}(t,x)=\varepsilon^{-3}
\zeta(t/\varepsilon^{2},x/\varepsilon)$. Next, for
locally summable $u(t,x)$ use the notation
$$
u ^{(\varepsilon)}(t,x)=u*\zeta_{\varepsilon}.
$$
The following properties of convolutions are well known
and are obtained by elementary manipulations.

\begin{lemma}
                               \label{lemma 6.25.1}
Assume that for   $x,y\in\bR$
and $s,t\in [0,1]$ we have
$$
|u(t,x)-u(t,y)|\leq |x-y|^{\beta} ,
$$
$$
|u(t,x)-u(s,x)|\leq |t-s|^{\beta/2}+a,
$$ 
 where $a\geq0$
is a fixed number.
Then for any $\varepsilon\in(0,1)$,
$x,y\in\bR$,
and $s,t\in [0,1-\varepsilon^{2}]$, we have
$$
|u ^{(\varepsilon)}(t,x)-u (t,x)|\leq  2\varepsilon^{\beta}+a,
$$

$$
|\partial^{2}_{t}u^{(\varepsilon)}(t,x)|+
|D^{4}u^{(\varepsilon)}(t,x)|+
|\partial_{t}D^{2}u^{(\varepsilon)}(t,x)|\leq
N\varepsilon^{-4}(\varepsilon^{\beta}+a),
$$

$$
|\partial_{t}u^{(\varepsilon)}(t,x)-\partial_{t}u^{(\varepsilon)}(s,x)|+
|D^{2}u^{(\varepsilon)}(t,x)-D^{2}u^{(\varepsilon)}(s,x)|\leq
N\varepsilon^{-2}[|t-s|^{\beta/2}+a],
$$

$$
|\partial_{t}u^{(\varepsilon)}(t,x)-\partial_{t}u^{(\varepsilon)}(t,y)|+
|D^{2}u^{(\varepsilon)}(t,x)-D^{2}u^{(\varepsilon)}(t,y)|\leq
N\varepsilon^{-2}|x-y|^{\beta},
$$
where $N$ is a constant depending only on $\zeta$.

\end{lemma}

\mysection{Proof of Theorem \protect\ref{theorem 6.22.1}}
                                            \label{section 6.26.1}

{\em Step 1. Estimating $v_{n}-v$ from above\/}.
By Lemma \ref{lemma 6.22.1} (a), if $\underline{\sigma}
\leq\sigma\leq \bar \sigma$ and 
 $t,\delta\in[0,1]$ are such that $t+\delta\leq1$, then
$$
v(t,x)\geq Ev(t+\delta, x+  \sigma w_{\delta}).
$$

It follows that,  if $t+\delta\leq 1-\varepsilon^{2}$, then
$$
v ^{(\varepsilon)}(t,x)\geq 
Ev ^{(\varepsilon)}(t+\delta, x+ \sigma w_{\delta}).
$$
Now, the arbitrariness of $\delta$ and the smoothness of $v ^{(\varepsilon)}$
imply that
\begin{equation}
                                                   \label{6.26.3}
\partial_{t}v ^{(\varepsilon)}+(1/2)\sigma^{2}
D^{2}v ^{(\varepsilon)}\leq0
\end{equation}
in $[0,1-\varepsilon^{2}]\times\bR$
as long as $\underline{\sigma}
\leq\sigma\leq \bar \sigma$.

We now claim that if $t\in T_{n}\cap(0,1-\varepsilon^{2}]$
and 
$$
a:=\sup_{x\in\bR}[v_{n}(t,x)-v ^{(\varepsilon)}(t,x)],
$$
then for any $s\leq t, s\in T_{n}$ we have
\begin{equation}
                                                   \label{6.25.3}
v_{n}(s,x)\leq v ^{(\varepsilon)} (s,y)+a+ M_{n,\varepsilon}(t-s),
\end{equation}
where $M_{n,\varepsilon}= 
 N_{1}\varepsilon^{-2}n^{-\beta/2} $
and the constant $N_{1}$ is specified below.

If $s=t$, \eqref{6.25.3} is given. If it holds for 
$0<s\leq t, s\in T_{n}$, then  
$$
v_{n}(s-1/n,x)=\sup_{\xi\in\Theta}E v_{n} (s,x+\xi/\sqrt{n})
$$
$$
\leq \sup_{\xi\in\Theta}E v^{(\varepsilon)} (s,x+\xi/\sqrt{n})
 +a+M_{n,\varepsilon}(t-s)
$$
$$
=\sup_{\xi\in\Theta}E \big[v^{(\varepsilon)}(s,x+\xi/\sqrt{n})
-v^{(\varepsilon)}(s-1/n,x)\big]
$$
\begin{equation}
                                               \label{6.26.2}
+v^{(\varepsilon)} (s-1/n,x)+a
+M_{n,\varepsilon}(t-s).
\end{equation}
Here by Taylor's formula
$$
v^{(\varepsilon)}(s,x+\xi/\sqrt{n})-v^{(\varepsilon)}(s-1/n,x)
=v^{(\varepsilon)}(s,x)-v^{(\varepsilon)}(s-1/n,x)
$$
\begin{equation}
                                               \label{6.28.2}
+
(\xi/\sqrt{n})Dv^{(\varepsilon)}(s,x)
+(1/2)(\xi^{2}/n)D^{2}v^{(\varepsilon)}(s,x)+I,
\end{equation}
where
$$
I:=(1/2)(\xi^{2}/n)[D^{2}v^{(\varepsilon)}(s,x
+\theta_{1}\xi/\sqrt{n})-D^{2}v^{(\varepsilon)}(s,x)], \quad \theta_{1}\in(0,1).
$$

By Lemma \ref{lemma 6.25.1}
$$
E|I|\leq N\varepsilon^{-2}n^{-1-\beta/2}E|\xi|^{2+\beta}
\leq N\varepsilon^{-2}n^{-1-\beta/2}.
$$

It follows that the last supremum in \eqref{6.26.2} is dominated by
$$
J+
(1/(2n))\sup_{\underline{\sigma}\leq\sigma\leq\bar\sigma}
\sigma^{2}D^{2}v^{(\varepsilon)}(s,x)
+N\varepsilon^{-2}n^{-1-\beta/2},
$$
where
$$
J:=v^{(\varepsilon)}(s,x)-v^{(\varepsilon)}(s-1/n,x)=(1/n)
\partial_{t}v^{(\varepsilon)}(s-\theta_{2}/n,x)
$$
with $\theta_{2}\in(0,1)$, so that by Lemma \ref{lemma 6.25.1}
$$
J\leq (1/2)\partial_{t}v^{(\varepsilon)}(s ,x)
+N\varepsilon^{-2}n^{-1-\beta/2}.
$$

Hence, from   \eqref{6.26.3} and \eqref{6.26.2} we infer
$$
v_{n}(s-1/n,x)\leq v(s-1/n,x)+a +N_{1}\varepsilon^{-2}n^{-1-\beta/2}
+M_{n,\varepsilon}(t-s)
$$
$$
=v(s-1/n,x)+a 
+M_{n,\varepsilon}(t-s+1/n).
$$
This proves \eqref{6.25.3}.

By taking $t$ as the largest element in $T_{n}$ which 
is $\leq
1-\varepsilon^{2}$ and observing that
$$
|v ^{(\varepsilon)}(t,x)-v(t,x)|\leq N\varepsilon^{\beta},
\quad |v  (t,x)-g(x)|\leq\bar{\sigma}^{\beta}|1-t|^{\beta/2}
\leq N(\varepsilon^{\beta}+n^{-\beta/2}),
$$
$$
|v_{n}(t,x)-g(x)|\leq\bar{\sigma}^{\beta}|1-t|^{\beta/2} \leq N(\varepsilon^{\beta}+n^{-\beta/2}),
$$
we conclude   from \eqref{6.25.3} that
$$
v_{n}(0,0)\leq v(0,0) 
+N\varepsilon^{-2}n^{-\beta/2}+N(\varepsilon^{\beta}+n^{-\beta/2}).
$$
By taking here $\varepsilon=n^{-\beta/(4+2\beta)}$ we get
$$
v_{n}(0,0)\leq v(0,0) 
+N n^{^{-\beta^{2}/(4+2\beta)}} .
$$

{\em Step 2. Estimating $v-v_{n}$ from above\/}.
It is convenient to extend 
  $v_{n}(t,x)$ defined on $T_{n}$ to the whole
of $[0,1]$ keeping its values on $T_{n}$ and making it
constant on each interval $[k/n,(k+1)/n)$ and equal there 
to $v_{n}(k/n,x)$. We keep the notation $v_{n}$ for the extended function.
Observe that for $x,y\in \bR$ and $s,t\in [0,1]$ we have
$$
|v_{n}(t,x)-v_{n}(s,x)|\leq |t-s|^{\beta/2}+n^{-\beta/2},
$$
\begin{equation}
                                   \label{6.26.6}
|v_{n}(t,x)-v_{n}(t,y)|\leq |x-y|^{\beta}.
\end{equation}

Then for 
 $t \in[0,1]$  such that $t+1/n\leq1$
and $\xi\in\Theta$, we have
$$
v_{n}(t,x)\geq Ev_{n}(t+1/n, x+  \xi/\sqrt{n}).
$$

It follows that  if $t+1/n\leq 1-\varepsilon^{2}$, then
$$
v_{n} ^{(\varepsilon)}(t,x)\geq 
Ev_{n} ^{(\varepsilon)}(t+1/n, x+ \xi/\sqrt{n}).
$$
 
Here
$$
v_{n} ^{(\varepsilon)}(t+1/n, x+ \xi/\sqrt{n})=
v_{n} ^{(\varepsilon)}(t+1/n, x)
$$
\begin{equation}
                                             \label{6.28.4}
+
(\xi/\sqrt{n})Dv_{n} ^{(\varepsilon)}(t+1/n, x)
+(\xi^{2}/(2n))D^{2}v_{n} ^{(\varepsilon)}(t+1/n, x)+I,
\end{equation}
where
$$
I:=(1/2)(\xi^{2}/n)[D^{2}v_{n}^{(\varepsilon)}(t+1/n,x
+\theta_{1}\xi/\sqrt{n})-D^{2}v_{n}^{(\varepsilon)}(t+1/n,x)], \quad \theta_{1}\in(0,1).
$$

Owing to \eqref{6.26.6} and Lemma \ref{lemma 6.25.1}
$$
E|I|\leq Nn^{-1-\beta/2}\varepsilon^{-2}.
$$
Also
$$
|D^{2}v_{n} ^{(\varepsilon)}(t+1/n, x)-
D^{2}v_{n} ^{(\varepsilon)}(t , x)|\leq N\varepsilon^{-2}n^{-\beta/2}.
$$
Hence,
$$
v_{n} ^{(\varepsilon)}(t,x)\geq v_{n} ^{(\varepsilon)}(t+1/n, x)
+(1/(2n))\sup_{\underline{\sigma}\leq\sigma\leq\bar\sigma}
\sigma^{2}D^{2}v_{n} ^{(\varepsilon)}(t , x)-Nn^{-1-\beta/2}\varepsilon^{-2},
$$
\begin{equation}
                                                 \label{6.26.7}
Nn^{-1-\beta/2}\varepsilon^{-2}\geq
(1/(2n))\sup_{\underline{\sigma}\leq\sigma\leq\bar\sigma}
\sigma^{2}D^{2}v_{n} ^{(\varepsilon)}(t , x) +J,
\end{equation}
where
$$
J:=v_{n} ^{(\varepsilon)}(t+1/n, x)-v_{n} ^{(\varepsilon)}(t,x)
=(1/n)\partial_{t}v_{n} ^{(\varepsilon)}(t,x)
$$
$$
+(1/n)[
\partial_{t}v_{n} ^{(\varepsilon)}(t+\theta/n, x)-\partial_{t}v_{n} ^{(\varepsilon)}(t,x)]
$$
with $\theta\in(0,1)$. It follows that
$$
J\geq (1/n)\partial_{t}v_{n} ^{(\varepsilon)}(t,x)
-Nn^{-1-\beta/2}\varepsilon^{-2},
$$
which along with \eqref{6.26.7} leads to
$$
\sup_{\underline{\sigma}\leq\sigma\leq\bar\sigma}
\sigma^{2}D^{2}v_{n} ^{(\varepsilon)}(t , x)+
\partial_{t}v_{n} ^{(\varepsilon)}(t,x)\leq Nn^{-\beta/2}\varepsilon^{-2}
$$
as long as $t \leq 1-\varepsilon^{2}-1/n=:S$.

Now by It\^o's formula for any $\sigma_{\!\cdot}\in\Sigma$
$$
Ev_{n} ^{(\varepsilon)}(S,x^{\sigma_{\!\cdot}}_{S}(0,0))
\leq v_{n} ^{(\varepsilon)}(0,0)
+Nn^{-\beta/2}\varepsilon^{-2}.
$$
Here
$$
|v_{n} ^{(\varepsilon)}(0,0)-v_{n} (0,0)|+
|v_{n} ^{(\varepsilon)}(S,x)-v_{n}(S,x)|
\leq N(\varepsilon^{\beta}+n^{-\beta/2}),
$$
$$
 |v_{n}(S,x)-g(x)|\leq N(\varepsilon^{\beta}+n^{-\beta/2}),
$$
$$
E\big|g(x^{\sigma_{\!\cdot}}_{S}(0,0))-g(x^{\sigma_{\!\cdot}}_{1}(0,0))|
\leq E|\int_{S}^{1}\sigma_{r}\,dw_{r}|^{\beta}
\leq N(1-S)^{\beta/2}\leq N(\varepsilon^{\beta}+n^{-\beta/2}).
$$

We conclude that
$$
Eg( x^{\sigma_{\!\cdot}}_{1}(0,0))
\leq v_{n}  (0,0)
+Nn^{-\beta/2}\varepsilon^{-2}+N(\varepsilon^{\beta}+n^{-\beta/2}).
$$

Since this is true for any $\sigma_{\!\cdot}\in\Sigma$,
$$
v(0,0)\leq v_{n}  (0,0)
+Nn^{-\beta/2}\varepsilon^{-2}+N(\varepsilon^{\beta}+n^{-\beta/2}),
$$
which for $\varepsilon=n^{-\beta/(4+2\beta)}$ yields
$$
v(0,0)\leq v_{n}(0,0) 
+N n^ {-\beta^{2}/(4+2\beta)}  .
$$

This proves the theorem.

\mysection{Better rates of convergence}
                               \label{section 6.28.1}
In addition to our basic assumptions from 
Section \ref{section 6.28.2} suppose that
\eqref{6.27.1} holds with $\beta=2$, \eqref{6.22.3}
holds with $\beta=1$, and
$$
E(\xi^{3})=0\quad\forall\xi\in\Theta.
$$
\begin{theorem}
                                       \label{theorem 6.28.1}
Under the above assumptions
$$
|v(0,0)-v_{n}(0,0)|
\leq Nn^ {-1/4} ,
$$
where $N$ is a constant depending only on $\bar\sigma$ and
$M_{4}$.
 
\end{theorem}

Proof. We follow the proof of Theorem 
\ref{theorem 6.22.1} and while estimating
$v_{n}-v$ from above we represent the remainder $I$
in the Taylor  formula \eqref{6.28.2} as
$$
I=(1/6) (\xi^{3}/n^{3/2})D^{3}v^{(\varepsilon)}(s,x)+
(1/24) (\xi^{4}/n^{2})D^{4}v^{(\varepsilon)}(s,x
+\theta_{1}\xi/\sqrt{n}),  
$$
where $\theta_{1}\in(0,1)$.
By assumption and Lemma \ref{lemma 6.25.1}
$$
EI\leq N\varepsilon^{-3}n^{-2}.
$$

Furthermore, we represent $J$ as
$$
J=(1/n)\partial_{t}v^{(\varepsilon)}(s,x)
+(1/(2n^{2}))\partial^{2}_{t}v^{(\varepsilon)}(s-\theta_{2}/n,x)
$$
with $\theta_{2}\in(0,1)$ so that by Lemma \ref{lemma 6.25.1}
$$
J\leq (1/n)\partial_{t}v^{(\varepsilon)}(s,x)
+N\varepsilon^{-3}n^{-2}.
$$

This yields, as in the first part of the proof of 
Theorem 
\ref{theorem 6.22.1}, that   
$$
v_{n}(0,0)\leq v(0,0) +N\varepsilon^{-3}n^{-1}
+N(\varepsilon +n^{-1/2}).
$$
By taking here $\varepsilon=n^{-1/4}$ we get
$$
v_{n}(0,0)\leq v(0,0) 
+N n^{^{-1/4}} .
$$

To prove the estimate $v-v_{n}$ from above
we follow the second part of the proof
of Theorem 
\ref{theorem 6.22.1} and again take more terms in
Taylor's formulas. Then the remainder $I$ in \eqref{6.28.4}
takes the form
$$
(1/6)(\xi^{3}/n^{3/2})D^{3}v_{n}^{(\varepsilon)}(t+1/n,x)
+(1/24)(\xi^{4}/n^{2})D^{4}v_{n}^{(\varepsilon)}(t+1/n,x
+\theta_{1}\xi/\sqrt{n}),
$$
where $\theta_{1}\in(0,1)$.
By assumption and Lemma \ref{lemma 6.25.1}
$$
EI\leq N\varepsilon^{-4}(\varepsilon+n^{-1/2})n^{-2}.
$$
Also
$$
|D^{2}v_{n} ^{(\varepsilon)}(t+1/n, x)-
D^{2}v_{n} ^{(\varepsilon)}(t , x)|
=(1/n)|\partial_{t}D^{2}v_{n} ^{(\varepsilon)}(t+\theta_{2}/n, x)|.
$$
where $\theta_{2}\in(0,1)$. By Lemma \ref{lemma 6.25.1}
$$
|D^{2}v_{n} ^{(\varepsilon)}(t+1/n, x)-
D^{2}v_{n} ^{(\varepsilon)}(t , x)|\leq
N\varepsilon^{-4}(\varepsilon+n^{-1/2})n^{-1}.
$$
This brings us to
\begin{equation}
                                                 \label{6.28.7}
 N\varepsilon^{-4}(\varepsilon+n^{-1/2})n^{-2}\geq
(1/(2n))\sup_{\underline{\sigma}\leq\sigma\leq\bar\sigma}
\sigma^{2}D^{2}v_{n} ^{(\varepsilon)}(t , x) +J,
\end{equation}
where
$$
J:=v_{n} ^{(\varepsilon)}(t+1/n, x)-v_{n} ^{(\varepsilon)}(t,x)
=(1/n)\partial_{t}v_{n} ^{(\varepsilon)}(t,x)
+(1/(2n^{2})\partial^{2}_{t}v_{n} ^{(\varepsilon)}(t+\theta_{3}/n, x)
$$
with $\theta_{3}\in(0,1)$. By Lemma \ref{lemma 6.25.1}
$$
J\geq (1/n)\partial_{t}v_{n} ^{(\varepsilon)}(t,x)
-N\varepsilon^{-4}(\varepsilon+n^{-1/2})n^{-2},
$$
which along with \eqref{6.28.7} leads to
$$
\sup_{\underline{\sigma}\leq\sigma\leq\bar\sigma}
\sigma^{2}D^{2}v_{n} ^{(\varepsilon)}(t , x)+
\partial_{t}v_{n} ^{(\varepsilon)}(t,x)\leq 
N\varepsilon^{-4}(\varepsilon+n^{-1/2})n^{-1},
$$
as long as $t \leq 1-\varepsilon^{2}-1/n=:S$.

 Then by repeating the end of the proof of Theorem 
\ref{theorem 6.22.1} we arrive at
$$
v(0,0)\leq v_{n}  (0,0)
+N\varepsilon^{-4}(\varepsilon+n^{-1/2})n^{-1}
+N(\varepsilon +n^{-1/2}),
$$
which for $\varepsilon=n^{-1/4}$ yields
$$
v(0,0)\leq v_{n}(0,0) 
+N n^{-1/4} .
$$
The theorem is proved.

{\em Conjecture\/}: The assertion of Theorem
\ref{theorem 6.28.1} is essentially sharp
in the sense that, if $\phi=|x|$ and $\Theta$
is a singleton $\{\xi\}$, where 
 $\xi$ is allowed to depend on $n$
and be  such that $P(\xi=\pm1)=n^{-1/2}$,
$P(\xi=0)=1-2n^{-1/2}$, which will reflect
in the notation as $v_{n}^{n}$ and $v^{n}$,
then 
$$
 n^{1/4}v^{n}(0,0)=E|w|,
\quad \lim_{n\to\infty}n^{1/4}v^{n}_{n}(0,0)=0,
$$
where $w$ is normal $(0,2)$ random variable.

\end{document}